\documentclass{svproc}
\usepackage{sidecap}
\usepackage{url}
\usepackage{import}
\usepackage{color}
\usepackage{amsmath}
\usepackage{graphicx}        

\begin{document}
\mainmatter              
\title{Towards a finite volume discretization of the atmospheric surface layer consistent with physical theory}
\titlerunning{FV discretization of the atmospheric surface layer}  
%
\author{Simon Cl\'ement\inst{1} \and Florian Lemari\'e\inst{1} \and Eric Blayo\inst{2}}
\authorrunning{Cl\'ement et al.} 
%
\tocauthor{}
\institute{
Univ. Grenoble Alpes, Inria, CNRS, Grenoble INP, LJK, 38000 Grenoble, France\\
\and
Univ. Grenoble Alpes, CNRS, Inria, Grenoble INP, LJK, 38000 Grenoble, France}

\maketitle              

\begin{abstract}
We study an atmospheric column and its discretization. Because of 
numerical considerations, the column must be divided 
into two parts: (1) a surface layer, excluded from the computational 
domain and parameterized, and (2) the rest of the column, 
which reacts more slowly to variations in surface conditions.
A usual practice in atmospheric models is to parameterize
the surface layer without excluding it from the computational domain,
leading to possible consistency issues.
We propose here to unify the two representations in a Finite Volume
discretization.
In order to do so, the reconstruction inside the first grid cell is performed 
using the particular functions involved in the parameterizations and not 
only with polynomials.
Using a consistency criterion, surface layer management strategies 
are compared in different physical situations.
%
%
\end{abstract}
\section{Introduction}
\label{sec:intro}
\begin{SCfigure}
	\centering
\begingroup%
  \makeatletter%
  \providecommand\color[2][]{%
    \errmessage{(Inkscape) Color is used for the text in Inkscape, but the package 'color.sty' is not loaded}%
    \renewcommand\color[2][]{}%
  }%
  \providecommand\transparent[1]{%
    \errmessage{(Inkscape) Transparency is used (non-zero) for the text in Inkscape, but the package 'transparent.sty' is not loaded}%
    \renewcommand\transparent[1]{}%
  }%
  \providecommand\rotatebox[2]{#2}%
  \newcommand*\fsize{\dimexpr\f@size pt\relax}%
  \newcommand*\lineheight[1]{\fontsize{\fsize}{#1\fsize}\selectfont}%
  \ifx\svgwidth\undefined%
    \setlength{\unitlength}{207.86245025bp}%
    \ifx\svgscale\undefined%
      \relax%
    \else%
      \setlength{\unitlength}{\unitlength * \real{\svgscale}}%
    \fi%
  \else%
    \setlength{\unitlength}{\svgwidth}%
  \fi%
  \global\let\svgwidth\undefined%
  \global\let\svgscale\undefined%
  \makeatother%
  \begin{picture}(1,0.30121875)%
    \lineheight{1}%
    \setlength\tabcolsep{0pt}%
    \put(0,0){\includegraphics[width=\unitlength,page=1]{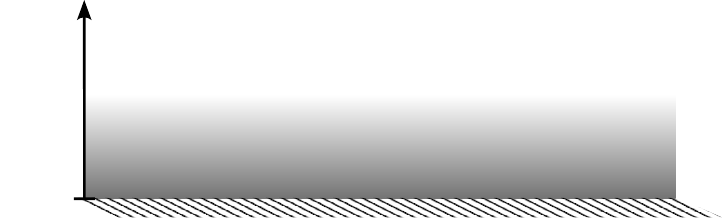}}%
    \put(-0.00374513,0.0217845){\color[rgb]{0,0,0}\makebox(0,0)[lt]{\lineheight{1.25}\smash{\begin{tabular}[t]{l}$0$\end{tabular}}}}%
    \put(-0.00374513,0.24549058){\color[rgb]{0,0,0}\makebox(0,0)[lt]{\lineheight{1.25}\smash{\begin{tabular}[t]{l}$z$\end{tabular}}}}%
    \put(0,0){\includegraphics[width=\unitlength,page=2]{notations_EkmanContinuous.pdf}}%
    \put(-0.00120603,0.12281306){\color[rgb]{0,0,0}\makebox(0,0)[lt]{\lineheight{1.25}\smash{\begin{tabular}[t]{l}$\delta_{\rm a}$\end{tabular}}}}%
    \put(0,0){\includegraphics[width=\unitlength,page=3]{notations_EkmanContinuous.pdf}}%
    \put(0.13171522,0.19061326){\color[rgb]{0,0,0}\makebox(0,0)[lt]{\lineheight{1.25}\smash{\begin{tabular}[t]{l}$(\partial_t + if) u - \partial_z (K_u \partial_z u) = if u_G$\end{tabular}}}}%
    \put(0.35584833,0.08545457){\color[rgb]{0,0,0}\makebox(0,0)[lt]{\lineheight{1.25}\smash{\begin{tabular}[t]{l}$K_u \partial_z u = u_\star^2 e_\tau$\end{tabular}}}}%
  \end{picture}%
\endgroup%

	\caption{Continuous equations in the computational domain  $(\delta_{\rm a}, +\infty)$ and constant flux in the SL $(0,\delta_{\rm a})$.
    The MO theory specifies the complex constant
    $u_\star e_\tau$ and a range of possible values
    for $\delta_{\rm a}$.
 }
	\label{fig:ND_NeutralCase_EkmanContinuous}
\end{SCfigure}

%
%
A common difficulty for atmospheric models is to represent the surface layer (SL), i.e. the area directly and almost instantaneously influenced by the presence of the ground or the ocean.
The scales in the SL (approximately the first 10 meters of the air column) are so small that the resolution needed for a numerical model to represent the phenomena correctly is out of reach. However, the Monin-Obukhov (MO) theory, which generalises the wall law to density-stratified fluids, provides under certain simple hypotheses (quasi-stationarity, horizontal homogeneity, etc.) an analytical formulation of the solution in the SL and the expression of the fluxes (heat, momentum) exchanged with the atmosphere above it. At the discrete level, the present treatment of this SL in numerical models is inconsistent: it is both treated like the rest of the atmosphere column by a standard numerical scheme (polynomial profile) when discretising the equations, and in a parameterised form (MO profile, which is a perturbation of a logarithmic profile,
see e.g. \cite{clement_phd}) for the calculation of fluxes. The consequences of this inconsistency are still poorly assessed in the context of atmospheric modeling, but we can mention for instance that in the context of combustion, it is mentioned in \cite{Jaegle_etal_2010} that the way the wall law 
is implemented in a given code and the way it interacts with the numerical methods 
used (in particular the turbulence scheme) can influence the numerical results as 
much as a particular choice of wall law. In this paper, we will address this inconsistency and propose a new finite volume formulation to remedy it.
\paragraph{The turbulent Ekman layer model}
Our approach is derived hereafter in the case of the 1D vertical Ekman layer model \cite{Mcwilliams_2009} in the neutral case. It includes the Coriolis effect with a constant parameter $f$ and a vertical turbulent flux term
$\langle w'u'\rangle$:
\begin{equation}
	\label{eq:ND_Intro_physics_evolutionEq}
	\partial_t u + i f u + \partial_z \langle w'u'\rangle
	= i f u_G
\end{equation}
where $i$ is the imaginary unit.
The constant nudging term $u_G$ pulls the solution towards the geostrophic equilibrium
(a large-scale solution where the pressure gradient balances the Coriolis force).
The horizontal wind $u$ (in ${\rm m}.{\rm s}^{-1}$) is a
complex variable accounting for both orientation and speed
of the wind.
The so-called Boussinesq hypothesis states that the turbulent flux is proportional to the gradient of $u$:
$\langle w'u' \rangle = - K_u \partial_z u$ where $K_u \geq 0$ is the turbulent viscosity. In the SL, the MO theory states that this turbulent flux is constant along the vertical axis and provides its analytical expression. We thus obtain  
the system of equations given in Fig.
\ref{fig:ND_NeutralCase_EkmanContinuous}. \par
\paragraph{Usual approach in atmospheric models.}
The usual choice made in current atmospheric models is to consider that the SL extends from the wall 
to the center of the first cell (let us note $z_{1/2}$ this altitude). In practice this 
means that compatibility constraints should apply at $z = z_{1/2}$ to correctly 
connect the profile as parameterized in the SL with the upper profile obtained using 
the numerical model. However the usual practice is to integrate the region 
$z \le z_{1/2}$ corresponding to the SL into the computational domain and then use 
the MO parameterization only to predict a surface flux at $z = 0$. The impact of this approximation is poorly 
documented to date: (i) the extent of the SL is fixed for purely numerical reasons and not for physical reasons, which does not guarantee that the solution will converge 
with the resolution; (ii) the solution in the area $z \le z_{1/2}$ is both parameterized and computed by the model, without ensuring the consistency of these two profiles. 
The coupling between the SL and the rest of the model is thus weak: the 
model provides the flow information at $z = z_{1/2}$ to the SL scheme and the 
latter provides in exchange a surface flux to the model at $z = 0$. In general 
no other interaction exists. For example, with this kind of coupling, the SL 
structures cannot really interact with the rest of the flow. 
For more details and numerous references, see \cite{Larsson_etal_2016}.\par 
A few studies address this issue: several alternatives for implementing a wall law in a Large Eddy Simulation solver 
are proposed in \cite{Jaegle_etal_2010};
%
%
a first step toward a proper Finite Volume (FV) approach is proposed by \cite{Nishizawa_Kitamura_2018}, where the authors extend the SL to the entire first cell $(0,z_1)$ and design a
scheme adapted to FV. \par 
In this paper, we propose to implement directly in the FV discretization the 
existing assumptions underlying the MO theory. In order to 
do so, the reconstruction inside the first grid cell is performed 
using the analytical functions involved in the wall laws and not 
only with polynomials. This approach also allows to relax the artificial assumption $\delta_a=z_{1/2}$ and to
extend the height of the SL beyond the first grid point if necessary (note that the log-layer mismatch, a well known numerical problem 
in Large Eddy Simulations, comes from a too thin SL).
By being able to choose the thickness of the SL based 
on physical --- and not only numerical --- criteria, the consistency of 
the schemes is improved. The vertical resolution can thus be refined 
without changing the continuous equations solved by the discretization, 
thus answering the issues raised in \cite{Basu_Lacser_2017}. 
Numerical experiments with a 1D Ekman layer model are performed, and SL  management strategies 
are compared for different types of stratification.
\section{The Finite Volume scheme}
\begin{SCfigure}
 \scalebox{0.9}{
	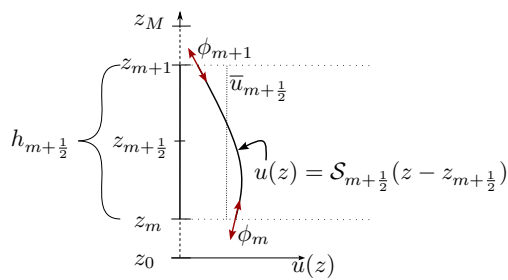}
	\caption{Summary of the notations related to the discretisation.}
	\label{fig:ND_NeutralCase_summary_notations}
\end{SCfigure}
\paragraph{Spline reconstruction of solutions.}
The space domain is divided into $M$ cells delimited by
heights $(z_0=0,\ldots, z_m,\ldots, z_M)$. The size of the $m$-th cell
is $h_{m+\frac{1}{2}}=z_{m+1}-z_{m}$ and the average of $u(z)$
over this cell is noted
$\overline{u}_{m+\frac{1}{2}}=\frac{1} {h_{m+\frac{1}{2}}}
\int_{z_{m}}^{z_{m+1}}u(z)dz$.
The space derivative of $u$ at $z_m$ is noted $\phi_{m}$.
Fig. \ref{fig:ND_NeutralCase_summary_notations} summarizes
these notations.
Averaging the evolution equation over a cell gives
the semi-discrete equation
\begin{equation}
\label{eq:ND_NeutralCase_semiDiscreteEkmanEq}
	(\partial_t + if) \overline{u}_{m+\frac{1}{2}} - 
	\frac{K_{u, m+1} \phi_{m+1} - K_{u, m} \phi_{m}}
		{h_{m+\frac{1}{2}}} = i f u_G
\end{equation}
The reconstruction of $u(z) = {\cal S}_{m+\frac{1}{2}}
				(z - z_{m+\frac{1}{2}})$
is chosen to be a quadratic polynomial (higher order schemes can be similarly derived, see
\cite{clement_phd}). The continuity
of $u(z)$ and its space derivative $\phi$ between cells yields the
relation:
\begin{equation}
\label{eq:ND_NeutralCase_continuityEquationFV}
	\frac{h_{m-1/2}}{{6}} \phi_{m-1} 
	+ \frac{h_{m+1/2}+h_{m-1/2}}
		{{3}} \phi_m  
	+ \frac{h_{m+1/2}}{{6}} \phi_{m+1} = \Bar{u}_{m+\frac{1}{2}} - \Bar{u}_{m-\frac{1}{2}}
\end{equation}
which is a FV approximation used in fourth-order
compact schemes for the first derivative $\partial_z u$ and
second-order for $\partial_z^2 u$ (e.g. \cite{piller_finite-volume_2004}).
\paragraph{Usual treatment of the SL with Finite Volume methods.}
The typical treatment of the SL in atmospheric models is to
use the evolution equation in the first cell $(z_0,z_1)$ to compute $\overline{u}_{\frac{1}{2}}$
and then assume that this averaged value is the wind speed at the center of the cell in
the Monin-Obukhov theory applied with $\delta_a=z_1$.
The corresponding bottom boundary condition is then
		  \begin{equation}
		\underbrace{K_{u,0} {\phi}^{n+1}_{0}}_{\text{Surface flux}}
		= u_\star^2 e_\tau \quad\hbox{with }
			u_\star = {\rm BULK}(
			\underbrace{\overline{u}^{n}_{\frac{1}{2}}}_{
				\text{Average around }
			z_{\frac{1}{2}}})
		  \end{equation}
where ${\rm BULK}$ is a routine based on the Monin-Obukhov theory, 
$e_{\tau} = \frac{\overline{u}_{\frac{1}{2}}^{n+1}}{||\overline{u}_{\frac{1}{2}}^{n}||}$, and $n$ denotes the time step.
This method has several drawbacks:
\begin{itemize}
    \item The value at the center of the cell is systematically
    larger than the average value because of the concavity of Monin-Obukhov
    profiles. This leads to a systematic underestimation of the surface flux by the SL scheme.
    A specific SL scheme was designed in
    \cite{Nishizawa_Kitamura_2018} to
    prevent this bias.
    \item The evolution equation is not compatible with the constant
    flux hypothesis that defines the SL, as introduced in Section \ref{sec:intro}.
    \item $\delta_{\rm a}$, the height of the SL,  is driven only
    by the space step and does not take into account any physical consideration.
\end{itemize}
\paragraph{On the incompatibility.}
\label{sec:ND_StratifiedCase_viscosity0_FVpure}
According to the wall law, $K_{u,0}$ should be equal to the (very small)
molecular viscosity $K_{mol} \approx {\rm 10^{-5}} \;{\rm m}^{2}.{\rm s}^{-1}$.
However, the boundary condition $K_{u,0} \phi_0 = u_\star^2 e_\tau$
does not have the same influence depending on the numerical scheme used to discretize \eqref{eq:ND_NeutralCase_semiDiscreteEkmanEq}:
\begin{itemize}
	\item \textbf{Finite Differences}:
Injecting the boundary condition in the evolution equation at the first grid level gives
		\begin{equation}
			(\partial_t+if) u_{1/2} = \frac{1}{h_{1/2}}
			\left(K_{u,1}\frac{u_{3/2} - u_{1/2}}{h_1}
			 - {u_\star^2 e_\tau} \right)
		\end{equation}
where one can see that the value $K_{u,0}$ does not intervene in the equation.
\item \textbf{Finite Volumes}:
applying $(\partial_t+if)$ to \eqref{eq:ND_NeutralCase_continuityEquationFV} and using the polynomial reconstruction and the equations of Fig. \ref{fig:ND_NeutralCase_EkmanContinuous}, one can see that the FV scheme implicitly uses
\begin{equation}
	(\partial_t+if) u(z_1) =
	\frac{K_{u,1} \phi_1 - {u_\star^2 e_\tau}}{h_{1/2}}
	+ (\partial_t+if)\left(\frac{\phi_1}{3}
	+ {\frac{u_\star^2 e_\tau}{6 K_{u,0}}}\right)h_{1/2}
\end{equation}
The (small) value of $K_{u,0}$ directly appears when we assume the
parabolic profile inside the first grid cell.
As a result, $u(z_1)$ scales with $\frac{1}{K_0}$ and
exhibits unreasonable values.
To obtain physically plausible profiles,
one can replace $K_{u,0}$ by $K_{u,\delta}$:
the wall law is then denied and
$(\partial_z u)(z_0)$ is multiplied
by $\frac{K_{mol}}{K_{u,\delta}}$.
Note that this problem would not occur 
if the simple FV approximation
$h_m \phi_m \approx \overline{u}_{m+\frac{1}{2}} - \overline{u}_{m-\frac{1}{2}}$ was used instead
of \eqref{eq:ND_NeutralCase_continuityEquationFV}.
\end{itemize}
\paragraph{Toward a Finite Volume scheme coherent with the physical theory.}
To address the drawbacks of the usual method presented above,
we now construct a numerical boundary condition that is coherent
with the continuous model with a free value of $\delta_{\rm a}$, named ``FV free":
\begin{equation}
	\label{eq:ND_NeutralCase_bottomCondFVFree}
	\underbrace{K_{u,\delta}\,
	{\phi}^{n+1}_{\delta}}_{\text{Flux at }
	\delta_{\rm a}}
		= u_\star^2\, e_\tau^{\rm free} \quad\hbox{with }
			  u_\star = {\rm BULK}(
			\underbrace{u^n(\delta_{\rm a})}_{
				\text{Reconstruction at }
			\delta_{\rm a}})
\end{equation}
where $e_\tau^{\rm free} =
\frac{u^{n+1}(\delta_{\rm a})}{||u^n(\delta_{\rm a})||}$
is the orientation of $u(\delta_{\rm a})$ obtained with the spline
reconstruction.
For the sake of simplicity, we assume in the following that $\delta_{\rm a} < z_1$
(this hypothesis being easily relaxed by using the Monin-Obukhov profiles
as the reconstruction in cells entirely contained in the SL).
In the first grid cell, we assume that
the constant flux hypothesis applies for $z<\delta_{\rm a}$
and we separate this cell into two parts: the surface layer $(0,\delta_{\rm a})$ and the ``sub-cell" $(\delta_{\rm a}, z_1)$.
This split corresponds to the change of governing equations
in Fig. \ref{fig:ND_NeutralCase_EkmanContinuous}.
Let $\widetilde{h} = z_1 - \delta_{\rm a}$  be the size of the
upper sub-cell $(\delta_{\rm a}, z_1)$
and $\widetilde{u} = \frac{1}{\widetilde{h}}\int_{\delta_{\rm a}}^{z_1}
u(z)dz$ be the corresponding averaged value of $u$.
The following subgrid reconstruction is used:
\begin{equation}
	\label{eq:ND_NeutralCase_fullReconstruction}
	u(z) = \begin{cases} \displaystyle
		{{\cal S}_{1/2}}\left(z - \frac{z_1 +
				\delta_{\rm a}}{2}\right),
		~~ z \geq \delta_{\rm a}
		\\[3mm] 
  \displaystyle
  \int_0^z
  \frac{u_\star^2\, e_\tau^{\rm free}}
  {K_{u,z'}}\,
  {\rm d}z',
		~~~~~~~~~~ z < \delta_{\rm a}
	\end{cases}
\end{equation}
where a closed-form of the integral
for $z < \delta_{\rm a}$ is given by MO theory.
The quadratic spline ${{\cal S}_{1/2}}$ used for reconstruction
is computed with the averaged value $\widetilde{u}$, the
size of the sub-cell $\widetilde{h}$ and the fluxes at the
extremities $\phi_{\delta}$ and $\phi_1$:
its definition 
${{\cal S}_{1/2}}(\xi) = 
{\widetilde{\color{black} u}} +
\frac{\phi_{1} + \phi_{\delta}}{2} \xi
+ \frac{\phi_{1} -
\phi_{\delta}}{2{\widetilde{\color{black} h}}}
\left(\xi^2 - \frac{{\widetilde{\color{black} h}}^2}{12}\right)$
is thus similar to the one in the other cells.

\section{Numerical Experiments}
The strategies to handle the SL are now compared through a test
of consistency: for several strategies,
the differences between a low-resolution and
a high-resolution simulations are compared. The smaller the difference
between the low-resolution and the high-resolution simulations, the better the consistency of the scheme.
The proposed strategy ``FV free" is compared with ``FV1" (the
typical current practice with Finite Volumes), ``FV2" (an intermediate
between ``FV free" and ``FV1": similar to ``FV free" but where the height of the surface layer $\delta_{\rm a}$ is set to $z_1$) and ``FD" (a Finite Difference reference).
\par
\begin{itemize}
    \item The turbulent viscosity is parameterized with a one-equation
    turbulence closure based on turbulent kinetic energy. The code is available
    at \cite{code_phd} and an in-depth description in \cite{clement_phd}. An Euler implicit time
    scheme integrates the model over a full day of simulation.
\item Parameters are $\Delta t = {\rm 30} \; {\rm s}$, $u_G = {\rm 8} \; {\rm m}.{\rm s}^{-1}$,
$f=10^{-4} \; {\rm s}^{-1}$
\item For ``FV free" the same $\delta_{\rm a}$ is used in both low-
and high-resolution simulations, whereas the resolution
imposes $\delta_{\rm a}$ in the other configurations.
\item The vertical levels of the low resolution simulation are taken as the 25 first of the 137-level configuration of the atmospheric model \textit{Integrated Forecasting System} at ECMWF (European Centre for Medium-Range Weather Forecasts). The high-resolution simulation has 3 times more cells: two grid levels are added in each of the low-resolution grid cells.
The usual SL strategies are designed for low-resolution configurations: the latter can hence be
considered as reference solutions, compared through
the sensitivity to the resolution.
\end{itemize}
\textbf{Neutral case: }
\begin{figure}
    \centering
    \includegraphics[scale=0.62]{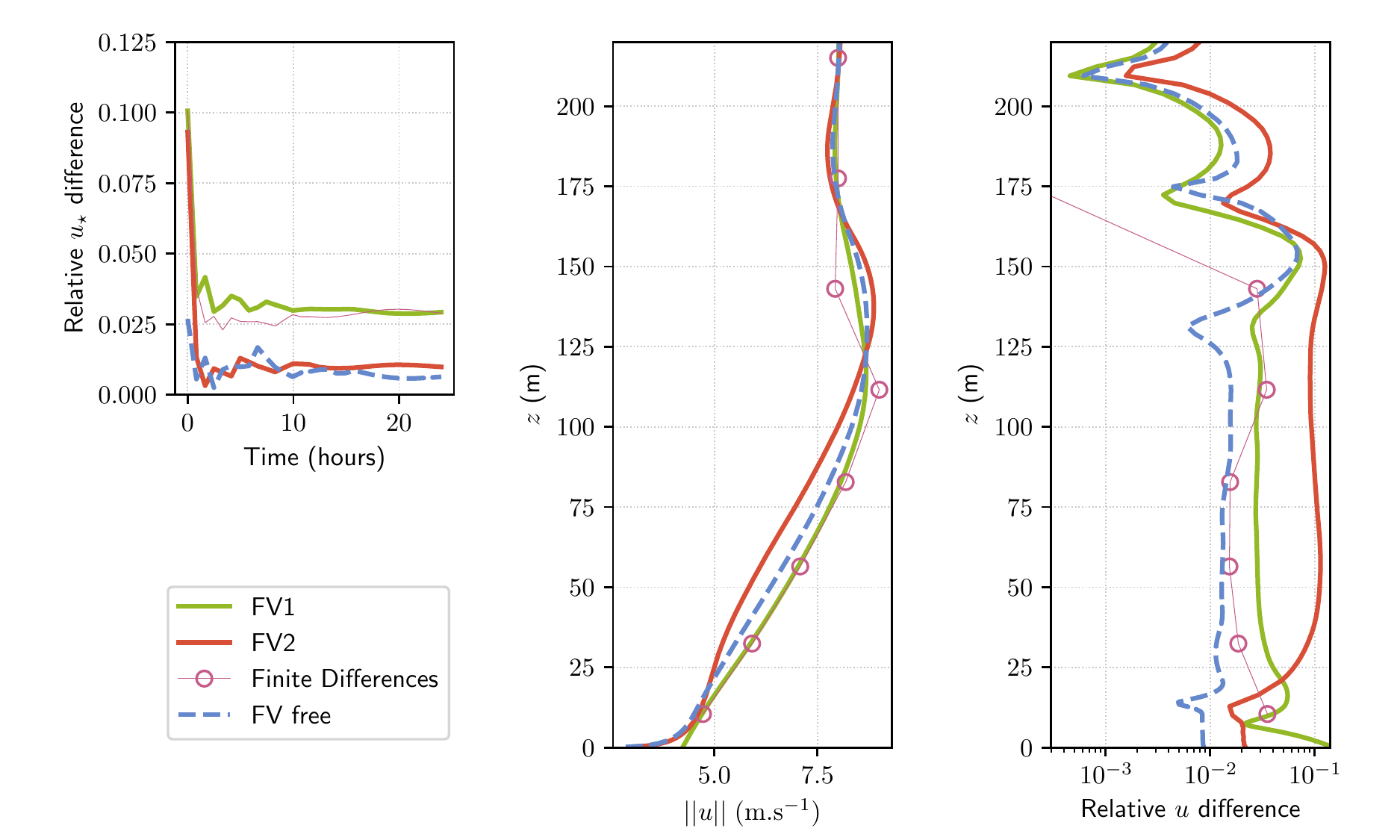}
    \caption{Relative difference between low-resolution and high-resolution
    simulations with several strategies for the SL handling. Left: Relative difference in $u_\star$ as a function of time. Center:
    vertical profiles of the wind speed at the end of the simulation. Right: Relative
    difference of the wind speed between low- and high-resolution
    along the vertical (note the log scale).
    }
\label{fig:consistency_comparisonNeutral}
\end{figure}
In the neutral case (constant density), 
the difference between low and high resolution of
the ``FV free" scheme is small at low altitude (see 
Fig. \ref{fig:consistency_comparisonNeutral}). This is mainly due to two factors:
\begin{itemize}
	\item $\delta_{\rm a} = z_{\frac{1}{2}}^{\rm low-res}$ is the same
		for both low and high resolutions, whereas for
		the other surface flux schemes the continuous
		equations change with $\delta_{\rm a}$.
	\item The initial relative difference for
		$u_\star$ (Fig. \ref{fig:consistency_comparisonNeutral}, left panel) is already much smaller than with the
		other schemes. This is a consequence of the
		imposed wall law: at initialization, there
		is already a logarithmic profile in the
		surface layer, instead of evolving toward a kind of compromise
  between the parameterized and the modeled values.
\end{itemize}

\textbf{Stratified case: }
We now focus on a stratified model \cite{Mcwilliams_2009} that includes
more of the physical behavior of atmospheric models:
the turbulence closure depends on the density $\rho$ such that
$\partial_z \rho \propto - \partial_z \theta$ where
$\theta$ is the potential temperature. We designed two
cases:
\begin{enumerate}
    \item A stable stratification, obtained with an initial temperature increasing  with the altitude, and a  surface temperature decreasing with time.
The initial potential temperature is 265 K in the first
100 meters of the atmosphere and then gains 1 degree every
100 meters; the surface temperature starts at 265 K and loses
1 degree every ten hours.
The ``low resolution" uses 15 grid points
in the 400m column and the ``high resolution"
uses 45 grid points.
\item An unstable stratification, obtained with a surface temperature following
a daily oscillation between $279$ K and $281$ K, and initial profiles of temperature and wind
 set to constant values of  $280$ K and
$8 \; {\rm m}.{\rm s}^{-1}$ respectively.
The ``low resolution" is composed of 50 grid levels of $10$m
each; 15 additional stretched levels between $500$m and
$1080$m make sure that the upper boundary condition
is not involved in the results.
The ``high resolution" divides every space cell into
3 new space cells of equal sizes.
\end{enumerate}
\par
\begin{SCfigure}
    \centering
    \includegraphics[scale=0.64]{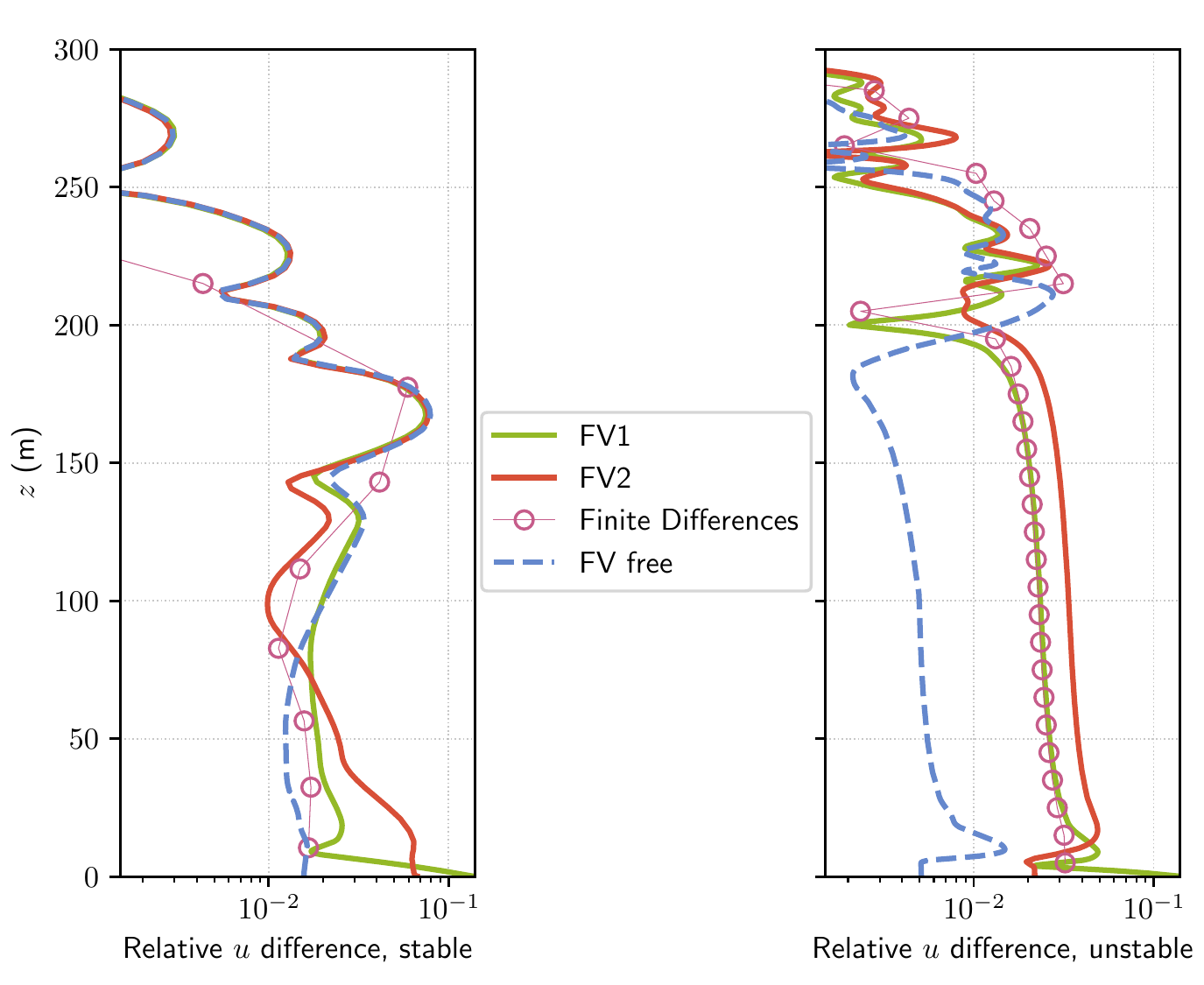}
    \caption{Relative difference of the wind speed between low- and high-resolution simulations for several SL strategies. left: stable stratification. Right: unstable stratification.}
    \label{fig:ND_Consistency_comparisonStratified}
\end{SCfigure}

 The differences between the two simulations are displayed in Fig.
\ref{fig:ND_Consistency_comparisonStratified}.
In the stable case,
the difference between the high resolution and the low resolution
results does not significantly change with the surface flux schemes.
The ``FV2" scheme is not very consistent because it tries to follow the
continuous model but with $\delta_{\rm a}$ changing  with the
resolution.
The Finite Difference or the ``FV1" methods suffer
less from this problem because, even if $\delta_{\rm a}$ changes,
it is assumed that the evolution equation is
integrated inside the surface layer. In
\cite{Maronga_etal_2020}, authors also find that the
sensitivity of their Large Eddy Simulation model to the grid spacing is
``\textit{more
likely related to under-resolved near-surface gradients
and turbulent mixing at the boundary-layer top, to the
[sub-grid scale] model formulation, and/or to numerical issues,
and not to deficiencies due to the use of improper surface
boundary conditions}".
\par
In the unstable case,  the  ``FV free" scheme
(with $\delta_{\rm a} = z_{\frac{1}{2}}^{\rm low-res}$) seems much
more robust than the other schemes in the first $200$ meters (remind that the height of the SL is approximately 10m).
Above this height the differences between the high
resolution and the low resolution simulations
are not clearly influenced by the SL treatment.
Note also that, as in the stable case, the ``FV2" scheme is also less consistent: enforcing the
MO theory in the first
cell increases the sensitivity of
the solution to $\delta_{\rm a}$ 
because the SL is then tightly coupled with the computational domain.
Finally, the ``FV free" scheme combines
good consistency properties with
a SL scheme coherent
with the physical theory.
%
%

\end{document}